\documentclass[a4paper, 11 pt]{article}
\usepackage{amsmath,amsthm,amsfonts,amssymb}%pacotes ams
\usepackage[mathcal]{eucal}
\usepackage[utf8]{inputenc}
\usepackage[T1]{fontenc}
\usepackage{graphicx}
\usepackage{wrapfig}
\usepackage{comment}
\usepackage{pdfpages}
\usepackage{indentfirst}
\usepackage[colorlinks,linkcolor=blue,anchorcolor=blue,,citecolor=blue,urlcolor=blue,backref,hyperindex]{hyperref}% Transforma referencias cruzadas em links.
\usepackage[top = 2 cm,left= 2 cm,right= 2 cm,bottom= 2 cm]{geometry}
\theoremstyle{plain}
\newtheorem{teo}{Teorema}

\theoremstyle{definition}

\theoremstyle{remark}

\title{Parametric Resonance of a charged pendulum with suspension point oscillating between two vertical charged lines}

\author{Carvalho, A. C.$^1$, Cabral, H. E.$^2$ and Araujo, G. C.$^3$}

\begin{document}
\maketitle

\centerline{$^{1}$ Adecarlos Costa Carvalho}
\medskip

{

       \centerline{Department of Mathematics}
   \centerline{Universidade Federal do Maranhão, Brazil}
     \centerline{E-mail: adecarlos.carvalho@ufma.br}
 } % Do not forget to end the \eightpoint by the sign } here
\bigskip

\centerline{$^{2}$ Hildeberto Eulalio Cabral}
\medskip

{

   \centerline{Department of Mathematics}
   \centerline{Universidade Federal de Pernambuco, Brazil}
      \centerline{E-mail: hild@dmat.ufpe.br}
 } % Do not forget to end the \eightpoint by the sign } here
\bigskip

 %% Enter the second author's name and address:

\centerline{$^{3}$ Gerson Cruz Araujo}
\medskip

{

   \centerline{Department of Mathematics}
   \centerline{Universidade Federal de Sergipe, Brazil}
   \centerline{S\~ao Cristov\~ao, Brazil}
   \centerline{E-mail:gerson@academico.ufs.br}
 } % Do not forget to end the \eightpoint by the sign } here
\bigskip

\begin{abstract}
In this work, we study a mathematical planar pendulum 
whose support point is
positioned 
equidistant between
two vertical and uniformly electrically charged wires. 
Its bob carries an electric charge and, its support point oscillates vertically, following a harmonic law of motion. 
%
%whose support point oscillates vertically, following an harmonic law of motion, and its bob carries an electric charge. The pendulum is positioned between two vertical wires carrying 
%
We study the dynamics of such phenomenon and the parametric resonances of the equilibria. Moreover, we obtain the surface in the parameter space (since such system presents three parameters) which separates the region of stability from the region of instability. On the particular case of zero charge, we obtain the boundary curves of the stability/instability of Matheiu equation.

\vspace{0.4 cm}


\vspace{0.4 cm}


\end{abstract}

\section{Introduction}

In Classical Mechanics, 
the mathematical planar pendulum
address to a model for a system consisting of 
a weight (the massive bob) suspended from a pivot (the support point) by a non flexible rod, so that, the bob can swing freely and the center of mass of the system is positioned at the bob.
Due to its huge number of variants and applications, this problem is known as one of the most studied problems in Mechanics. 
The non fixed support point variant of this phenomenon has received much attention, as we can see in \cite{Araujo_Cabral, Araujo_Cabral2, Bardin_Markeev, Cabral_Carvalho, Formalskii, Kholostova1, Kholostova2, Madigan, Neishtadt}.

In this paper, we approach the case where the support point, $O$, oscillates vertically in a harmonic way, the bob presents an electric charge $q$ and the system is positioned equidistant between two straight vertical wires uniformly electrically charged (see Fig. \ref{fig1}). The system, thus, presents three parameters: $\varepsilon$, a small parameter associated to the amplitude of the swing and the pendulum length, $\alpha$, a parameter associated to the oscillation frequency and the pendulum length, and last, but not least, $\mu$, a parameter associated to the pendulum length and the electric charges at the bob and at the vertical wires. For $\varepsilon = 0$, the support point is fixed.
In the case $\mu = 0$, the dynamics corresponds to a pendulum electrically charged whose support point is oscillating vertically between two vertical electrically charged wires, following a harmonic move.   

The problem we address in this work is described in its Hamiltonian formulation by (\ref{eq:hamiltonian1}). It presents two equilibria, $P_1 = (0, 0)$ and $P_2 = (\pi, 0)$, located at the vertical line containing the suspension point $O$. According to the value of parameter $\mu$, a equilibrium of the system may be stable or unstable (see Section \ref{S2}). In Section \ref{S3}, we normalize the linearized Hamiltonian at each equilibrium, restricted to the region of the parameter space where the equilibria are linearly stable. In Section \ref{S4}, we approach the boundary surfaces which separate the regions of stability and instability at the parameter space $(\mu, \alpha, \varepsilon)$. We calculate the coefficients of its parametrizations until the fifth order in terms of parameter $\mu$. 
By analyzing planar sections, $\mu = $ constant, we obtain the border curves of stability/instability in each plane. 
In the particular case $\mu = 0$, we obtain the border curves of Mathieu equation, matching the obtained coefficients with the ones found in \cite{Bardin_Markeev}.

%Em suma, o pêndulo com ponto de suspensão oscilante tem devida atenção desde os artigos advindos das pesquisas de Kapitsa \cite{Kapitza3, Kapitza4}. Seu estudo foi feito a partir de diferentes perspectivas e continua a atrair o interesse de pesquisadores de todo o mundo, consulte as referências \cite{Dadfar,Formalskii, Kallu, Kholostova1, Kholostova2, Levi, Morozov, Neishtadt, Ovseyevich}. Neste trabalho consideramos além da oscilação do ponto de suspensão, acrescentamos o efeito na dinâmica proveniente das cargas eletrostáticas, situadas no bulbo e nas linhas verticais eletrizadas uniformemente.

\section{Problem Formulation} \label{S2}

Consider a pendulum of length $l$, whose suspension point, $O$ is under a vertical harmonic oscillation described by the equation $\rho = a\cos \nu t$, $a > 0$. The pendulum bob has mass $m$ and is electrically charged. The pendulum swings between two vertical wires both uniformly charged with the same constant linear density. The pendulum support point remains equidistant between the two charged wires and we call $d$ this constant distance.

\begin{figure}[!ht]
\centering
\includegraphics[scale = 0.45]{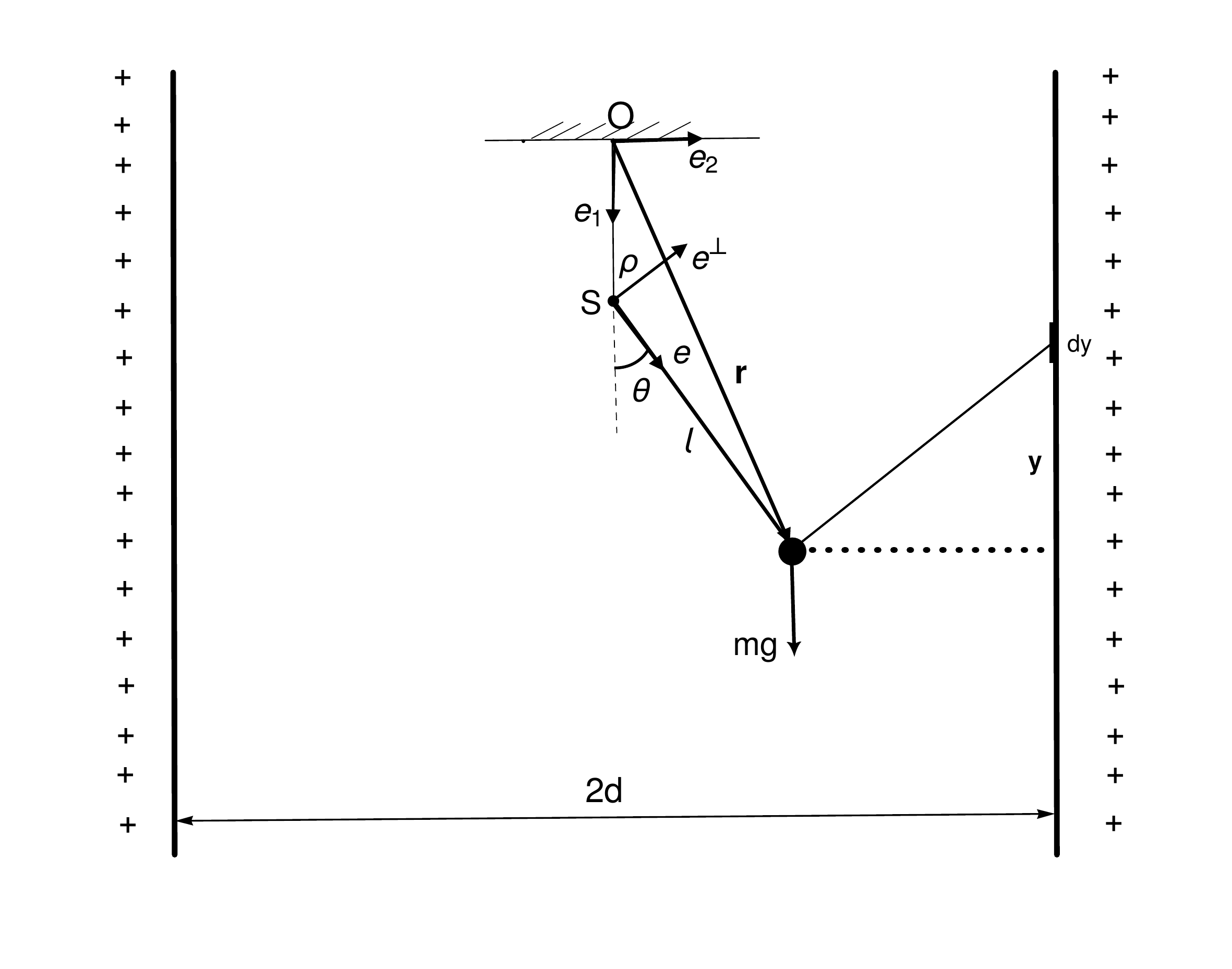}\caption{Charged pendulum with support point oscillating equidistant between two charged vertical wires}
\label{fig1}
\end{figure}

Let $y$ be the distance between an infinitesimal element $\mathrm{d}y$ to the orthogonal projection of the bob with respect to the same line and, $\theta$ be the angle in radians between the pendulum rod and the vertical direction (see Fig. \ref{fig1}). Both wires are considered having the same constant linear charge density $\sigma =\frac{\mathrm{d}Q}{\mathrm{d}y}$. Consider the inertial orthonormal basis $\mathbf{e}_1, \mathbf{e}_2$ at the suspension point $O$, where the first vector points down and the second vector points right. In this way, the position vector of the suspension point is given by $\textbf{R} = \rho \textbf{e}_1$. Consider $\textbf{r}$ the bob position vector and $(\textbf{e}, \textbf{e}^{\perp})$ a moving positive orthogonal basis satisfying $\textbf{e} = \frac{1}{l}(\textbf{r} - \textbf{R})$. 

The forces acting in the system are
\begin{center}
 $\textbf{F}_g = mg \textbf{e}_1$, \ \ $\textbf{T} = T\textbf{e}$,
\end{center}
\begin{eqnarray*}
\textbf{F}_r &=& \int_{-\infty}^{+\infty}\frac{-k_0 q \sigma dy}{y^2 + (d - l\sin \theta )^2}\cdot \frac{y \textbf{e}_1 + (d - l\sin\theta) \textbf{e}_2 }{\sqrt{y^2 + (d - l\sin \theta)^2}} =  -\frac{2 k_0 q  \sigma }{d - l\sin \theta} \textbf{e}_2,\\
\textbf{F}_l &=& \int_{-\infty}^{+\infty}\frac{- k_0 q \sigma dy}{y^2 + (d + l\sin \theta )^2}\cdot \frac{y \textbf{e}_1 - (d + l\sin\theta) \textbf{e}_2 }{\sqrt{y^2 + (d + l\sin \theta)^2}}= \frac{2 k_0 q  \sigma }{d + l\sin \theta}\textbf{e}_2,
\end{eqnarray*}
where $\textbf{F}_g$ is the gravitational force, $\textbf{T}$ is the tension at the rod connected to to the mass $m$, $\textbf{F}_l$ and $\textbf{F}_r$ are electrostatic forces due to the vertical charged wires from the left and from the right respectively and, $k_0$ is the Coulomb constant. We discard the influence of magnetic force since the velocity of the bob is considered to be low.   

Since $\textbf{e} = \cos \theta \textbf{e}_1 + \sin\theta \textbf{e}_2$ and $\textbf{e}^{\perp} = -\sin \theta \textbf{e}_1 + \cos\theta \textbf{e}_2$, we have $\dot{\textbf{e}} = \dot{\theta} \textbf{e}^{\perp}$ and $\textbf{e}^{\perp} = - \dot{\theta} \textbf{e}$. Moreover, we have, by taking the second derivative of, $\textbf{r} = \rho \textbf{e}_1 + l \textbf{e}$, 
\[\ddot{\textbf{r}} = \ddot{\rho}\textbf{e}_1 + l \ddot{\theta}\textbf{e}^{\perp} - l \dot{\theta}^2\textbf{e}.\]
The total force acting at the bob is given by $\textbf{F} = \textbf{F}_g + \textbf{T} + \textbf{F}_r + \textbf{F}_l$. Since $\textbf{e}_1 = \cos \theta \textbf{e} - \sin\theta \textbf{e}^{\perp}$, by comparing the components $\textbf{e}^{\perp}$ and applying the Newton's Second Law of Motion, $m\ddot{\textbf{r}} = \textbf{F}$, we obtain an expression for the pendulum movement. 
\[m l \ddot{\theta} -m\ddot{\rho}\sin\theta + mg \sin\theta + 2 k_0 q \sigma \left(\frac{ 1 }{d - l \sin\theta} - \frac{1}{d +l \sin\theta}\right)\cos\theta = 0 .\]

Let $\rho = a\cos \nu t$ be the harmonic law describing the suspension point movement, we found convenient to take $\rho = \nu \tau$ as a new independent variable. Thus, we can wright $\ddot{\theta} = \nu^2\theta''$ and $\ddot{\rho} = \nu^2 \rho''$, where the apostrophe notation indicates the derivative with respect to $\tau$. Now, by taking $d = 2l$,  $\alpha = \frac{g}{l\nu^2}$, $\varepsilon = \frac{a}{l}$, $\mu = \frac{q}{l\sigma}$  and $\sigma $ such that $\frac{4 k_0 \sigma^2}{m l \nu^2} = 1$, the movement equation becomes
\begin{equation}\label{eq-segunda-ordem-1}
	\theta''  + \alpha \sin\theta + \varepsilon \cos \tau \sin\theta + \frac{\mu  }{7 +  \cos 2\theta}\sin 2\theta = 0.
\end{equation}

%\begin{equation}
% \theta''  + \alpha \sin\theta + \varepsilon \cos \tau \sin\theta + \frac{ 2 \alpha q  }{7 +  \cos 2\theta}\sin 2\theta = 0 
%\end{equation}

For $\mu = 0$, equation \eqref{eq-segunda-ordem-1} describes the motion of a pendulum whose support point oscillates vertically in a harmonic way (see \cite{Bardin_Markeev}). For  $\theta = 0 $ and $\theta = \pi$ we have equilibrium points for every choice of the parameters. In this work, we present a study on the parametric resonance of these equilibria by means of the Hamiltonian formulation. 

%\section{Hamiltonian Formulation}

Writing $x = \theta$ and $y = \theta'$, we get a Hamiltonian System whose Hamiltonian function is given in terms of 
\begin{equation}\label{eq:hamiltonian1}
	 H(x,y,\tau, \mu, \alpha, \varepsilon) = \frac{1}{2}y^2 - \alpha \cos x - \varepsilon \cos \tau \cos x - \frac{\mu}{2} \log (7 + \cos 2 x). 
\end{equation}
The points $P_1=(0,0)$ and $P_2 = (\pi, 0)$ are equilibria of the system for an arbitrary choice of the parameters. For $\varepsilon = 0$, we have an autonomous system with one degree of freedom. For both equilibria we have $H_{yy} = 1$ and $H_{xy} = 0$.  The value of $H_{xx}$ is
\[
H_{xx} =  \frac{\mu}{4} + \alpha \quad \text{for } P_1 \qquad \text{ and } \qquad H_{xx} = \frac{\mu}{4} - \alpha \quad \text{for } \quad P_2.
\]
Thus, the equilibrium $P_1$ is stable whenever $\mu > -4 \alpha$ and unstable for $\mu < -4 \alpha$ while the equilibrium $P_2$ is stable when $\mu > 4\alpha$ and unstable when $\mu < 4 \alpha $.   

\section{Parametric Resonance} \label{S3}

This section is devoted to the study of the parametric resonance of the linearized system related to (\ref{eq:hamiltonian1}) at the parameter space $(\mu, \alpha, \varepsilon)$. We perform this study by applying the following result, know as Krein-Gelfand-Lidskii Theorem \cite{Markeev1}.

\begin{teo}[Krein-Gelfand-Lidskii] \label{theorem:krein}
	Given a linear Hamiltonian System whose Hamiltonian function is given by
	\begin{equation} \label{eq:hamiltonian_theorem_krein}
		H = \frac{1}{2} \sum_{k = 1}^n \sigma_k(x_k^2 + y_k^2 ) + \varepsilon H_1 + \varepsilon^2 H_2 + \cdots, 
	\end{equation}
	where $H_1, H_2, \cdots $ are quadratic forms with respect to $x_1, y_1, \cdots, x_n, y_n$; its coefficients are continuous and $2 \pi $ periodic in $t$. For $\varepsilon > 0 $ small enough, the linear system whose Hamiltonian is given by (\ref{eq:hamiltonian_theorem_krein}) is stable if, and only if, the terms $\sigma_k$ do not satisfy
	\begin{equation} \label{eq:theorem_krein}
		\sigma_k + \sigma_l = N,
	\end{equation}  
	for $k, l = 1, 2, \cdots, n$\ \ and \ \ $N = \pm 1, \pm 2, \cdots $.
\end{teo}

The system given in (\ref{eq:hamiltonian1}) is a time-depending Hamiltonian system.  It is also dependent of the parameters $\mu, \alpha$ and $\varepsilon$. On the following, we present a study on the stability of the linearized system around the equilibria $P_1 = (0,0)$ for $\mu > - 4\alpha$ and $P_2=(\pi, 0)$ for $\mu > 4\mu$.

Let us consider $\xi = x - x_0$ and $\eta = y$, with $x_0 = 0$ for $P_1$ and $x_0 = \pi$ for $P_2$ thus, the linearized Hamiltonian functions can be written as
\begin{eqnarray}
	H(\xi,\eta,\tau, \mu, \alpha, \varepsilon) &=& \frac{1}{2}\eta^2 + \frac{1}{2}\left[\varepsilon\cos\tau + \alpha + \frac{\mu}{4}  \right]\xi^2 \quad \mbox{for} \quad P_1, \label{eq:hamiltonian2_p1}\\
	H(\xi,\eta,\tau, \mu, \alpha, \varepsilon) &=& \frac{1}{2}\eta^2 - \frac{1}{2}\left[\varepsilon\cos\tau + \alpha - \frac{\mu}{4} \right]\xi^2 \quad \mbox{for} \quad P_2.\label{eq:hamiltonian2_p2} 
\end{eqnarray}
Now, by applying the symplectic change of coordinates $\xi, \eta \to x, y$ given by
$$\xi = \omega^{-1/4 } x ,\qquad  \eta = \omega^{1/4} y,$$
on (\ref{eq:hamiltonian2_p1}) and (\ref{eq:hamiltonian2_p2}) and expanding it into power series on $\varepsilon$, we obtain
\begin{equation} \label{eq:hamiltonian3}
	H(x,y,\tau, \mu, \alpha, \varepsilon) = \frac{\omega}{2}(x^2 + y^2) +  \frac{x^2 \cos \tau  }{2\omega}\varepsilon, 
\end{equation}
where, $\omega^2 = \frac{\mu}{4} + \alpha $ for $P_1$ and $\omega^2 = \frac{\mu}{4} - \alpha $ for $P_2$.

Note that the Hamiltonian (\ref{eq:hamiltonian3}) is on the form (\ref{eq:hamiltonian_theorem_krein}), where the frequency of the linear system, $\omega$, depends on the parameters $\mu$ and $ \alpha$. If, for some integer $N$, we have $2 \omega(\mu, \alpha) = N$ then, it follows from the Krein-Gelfand-Lidskii Theorem, that the linear unperturbed system is not stable. Moreover, by analyzing the system on the parameter space $(\mu, \alpha, \varepsilon)$, we see that the equation $2\omega(\mu, \alpha) = N$ defines a curve at the subspace $(\mu, \alpha, 0)$. Thus, for every $(\mu_0, \alpha_0, 0)$ in this curve if, $\varepsilon > 0$ then, $(\mu_0, \alpha_0, \varepsilon)$ may or may not be stable. 
In this way, we obtain boundary surfaces separating the regions of stability and instability in the parameter space $(\mu, \alpha, \varepsilon)$. Such surfaces will be expressed as a graph of a function under the plane $(\mu, 0, \varepsilon)$, as a power series on $\varepsilon$ and its coefficients will be given in terms of $\mu$, that is, 
\begin{equation}\label{eq:alpha}
	\alpha = \alpha_0 + \alpha_1 \varepsilon + \alpha_2 \varepsilon^2 + \alpha_3 \varepsilon^3 + \alpha_4 \varepsilon^4 + \mathcal{O}(\varepsilon^5),
\end{equation}
where $\alpha_{j} = \alpha_{j}(\mu)$, $j \geq 1$, are determined by $\mu$, $\alpha_0 = (N^2 -  \mu)/4$ for $P_1$ and $\alpha_0 = (\mu - N^2)/4$ for $P_2$ are curves from $(\mu,\alpha, 0)$ plane defined by the condition $2\omega(\mu, \alpha) = N.$ The coefficients $\alpha_j$, $j \geq 1$, are give in the next section.

We now apply (\ref{eq:alpha}) on the Hamiltonian terms (\ref{eq:hamiltonian2_p1}) and (\ref{eq:hamiltonian2_p2}) and, perform the symplectic change of coordinates given by 
$$\xi = \omega_0^{-1/4}\widetilde{X},  \quad \eta = \omega_0^{1/4}\widetilde{Y}$$
and
$$\widetilde{X} = X \cos (\sqrt{\omega_0}\tau)   + Y \sin (\sqrt{\omega_0}\tau), \quad \widetilde{Y} =  - X \sin (\sqrt{\omega_0}\tau)  + Y \cos (\sqrt{\omega_0}\tau) ,$$
in order to write 
 \begin{equation} \label{eq:hamiltonian4}
  \mathcal{H}(X,Y,\tau, \mu, \alpha, \varepsilon) = \frac{(-1)^{i+1}S^2 \cos\tau}{ N}\varepsilon + \sum_{j \geq 1}  \frac{\alpha_j S^2}{4 N} \varepsilon^j,
 \end{equation}  
where $i = 1$ for $P_1$, $i=2$ for $P_2$ and $S = X \cos (N\tau/2)   + Y \sin (N\tau/2)$. The rotation considered eliminates the term $\mathcal{H}_0$ from the Hamiltonian expression (see \cite{Cabral_Carvalho}), enabling a shorter computation of the coefficients at the boundary surfaces.

\section{Boundary surface of the stability/instability regions} \label{S4}

In this section, we apply the Depri-Hori Method \cite{Araujo_Cabral, Lucia_Hild, Kamel1, Markeev1} on the Hamiltonian function (\ref{eq:hamiltonian4}) in order to obtain boundary surfaces separating the stability and instability regions. This method enables us to transform, via simplectic change of variables $X,Y \to p, P$, Hamiltonian functions of the form
\begin{equation} \label{hamiltoniano-metodo-depri-hori}
	H(X, Y, \nu, \varepsilon) = \sum_{m = 0}^{\infty} \frac{\varepsilon^m}{m!} H_m (X, Y, \nu),
\end{equation}
into an autonomous Hamiltonian of the form
\begin{equation} \label{eq:hamiltonian_aut}
	K(p, P) = k_{02}p^2 + k_{11} p P + k_{20} P^2,
\end{equation}
where $k_{i j} = \sum_{m =1}^{\infty} k_{ij}^{(m) } \varepsilon^m$, with $k_{ij}^{(m)}$ depending on $\alpha_1, \cdots, \alpha_m$.

After applying the Depri-Hori Method to the Hamiltonian (\ref{eq:hamiltonian4}), we obtain a Hamiltonian of the form (\ref{eq:hamiltonian_aut}), whose term $k_{11 }$ is null. As a consequence, the characteristic equation is given by $\lambda^2 + 4 k_{20}k_{02} = 0.$ Thus, the stability region is determined by the condition $k_{20}k_{02} > 0$ and then, the boundary region is determined by the equation $k_{20}k_{02} = 0$, that is, 
\[
k_{20} = 0 \quad \text{ or } \quad k_{02} = 0.    
\]
The coefficients $\alpha_j(\mu)$ in \eqref{eq:alpha} can be determined by taking the coefficients of every power of $\varepsilon$ in the expressions $k_{20}$ and $k_{02}$ equals to zero and then, the boundary surfaces on the parameter space $(\mu, \alpha, \varepsilon)$.  The surfaces emanate from the curve $\alpha_0 = (N^2- \mu)/4$ for $P_1$ and $\alpha_0 = (\mu - N^2)/4 $ for $P_2$ given by the equation $2 \omega(\mu,\alpha) = N$, $N \ge 1$, at the plane $\varepsilon = 0$.

\subsection{Boundary surfaces for the equilibrium P1} 

In this subsection we present the surfaces which separates the regions of stability and instability for the equilibrium $P_1$. For each resonance of the form $2 \omega = N$, $N =1, 2, 3, \cdots$, by applying Depri-Hori Method to Hamiltonian (\ref{eq:hamiltonian4}), with $\omega_0 = \frac{\mu}{4} + \alpha_0$ and $N = 1$, we get
\begin{eqnarray*}
	k_{20}^{(1)} &=& \frac{1}{4} + \frac{\alpha_1}{2}, \\
	k_{20}^{(2)} &=& -\frac{1}{16}(3 + 12 \alpha_1 + 8 \alpha_1^2 - 8 \alpha_2),\\
	k_{20}^{(3)} &=& -\frac{3}{64} + \frac{\alpha_1}{2} + \frac{3 \alpha_1^2}{2} + \alpha_1^3 - \frac{3 \alpha_2}{4} - \alpha_1 \alpha_2 + \frac{\alpha_3}{2},
\end{eqnarray*}
and,
\begin{center}
	$k_{02}^{(1)} - k_{20}^{(1)} = -\frac{1}{2}$,\qquad $k_{02}^{(2)} - k_{20}^{(2)} = \frac{3 \alpha_1}{8}$,\\
	\vspace{0.5 cm}
	$k_{02}^{(3)} - k_{20}^{(3)} = \frac{3}{32} - 3 \alpha_1^2 + \frac{3 \alpha_2}{2}$.
\end{center}
The boundary surfaces on the parameter space are those for which the coefficients of $k_{20} $ and $k_{02}$ are zero. Thus, the boundary surfaces are the two given in parametric form by
\[
\alpha = \frac{1 - \mu}{4} \mp  \frac{1}{2}\varepsilon - \frac{1}{8} \varepsilon^2 \pm \frac{1}{32} \varepsilon^3 -\frac{1}{384} \varepsilon^4  \mp \frac{11}{4608}\varepsilon^5 + \mathcal{O}(\varepsilon^6).
\]
%\begin{figure}[ht]
%	\centering
	%  \includegraphics[scale = 0.8]{figuras/fig2v3.pdf}
%	\caption{Boundary surface for $P_1$ for $N =1$; Planar section for $\mu = -\frac{1}{2}.$}
%	\label{fig2}
%\end{figure}

 \begin{figure}[ht]
 \centering
 \includegraphics[scale = 0.8]{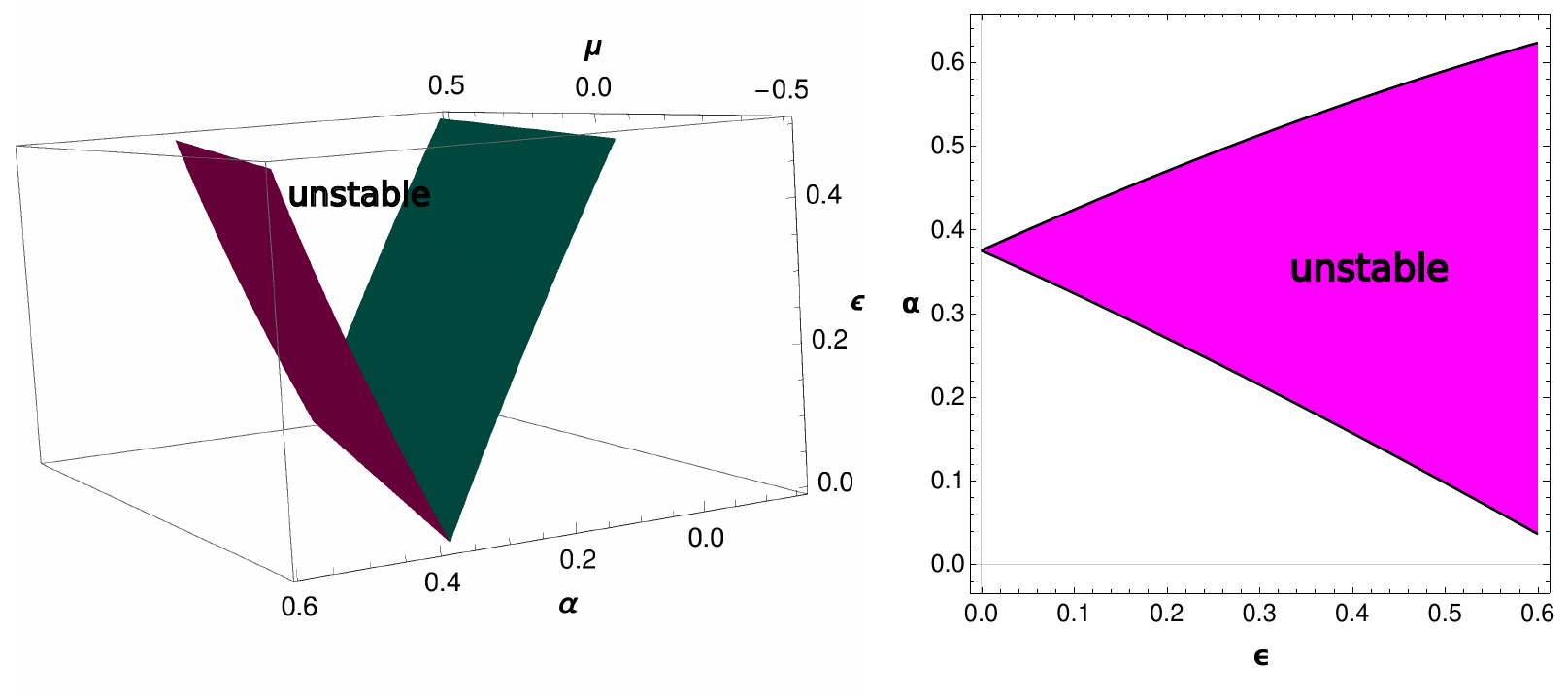}
  \caption{Boundary surface for $P_1$ and $N =1$; Planar section $\mu = -\frac{1}{2}.$}
  \label{fig2}
 \end{figure}

Both surfaces we found delimits the instability region related to the resonance $2 \omega = 1$ for the equilibrium  $P_1$. In Fig. \ref{fig2}, the picture on the left shows the surface delimiting the instability region. On the right hand side, we see the planar section of this surface whith respect to the plane $\mu = -\frac{1}{2}$.

On the resonance for $N = 2$, the equations $k_{20}$ and $k_{02}$ can be expressed as
\begin{eqnarray*}
 k_{20}^{(1)} &=& \frac{\alpha_1}{4},\\
 k_{20}^{(2)} &=& -\frac{1}{48}(5 + 9 \alpha_1^2 + 12 \alpha_2),\\
 k_{20}^{(3)} &=& \frac{1}{288}(70 \alpha_1 + 45 \alpha_1^3 -108 \alpha_1\alpha_2 + 72 \alpha_3),
\end{eqnarray*}
and,
\begin{center}
 $k_{02}^{(1)} = k_{20}^{(1)}$, \quad $k_{02}^{(2)} - k_{20}^{(2)} = \frac{1}{8} + \frac{\alpha_1^2}{4}$, \quad $k_{02}^{(3)} - k_{20}^{(3)} = -\frac{\alpha_1 }{48} (13 + 9 \alpha_1^2 - 24 \alpha_2)$.
\end{center}
In this case, the equations $k_{20} = 0$ and $k_{02} = 0$ provide, respectively the following surfaces
\begin{eqnarray*}
 \alpha &=& \frac{4 - \mu}{4} + \frac{5}{12}\varepsilon^2 - \frac{763}{3456 } \varepsilon^4+ \frac{1 002 401}{4 976 640} \varepsilon^6 + \mathcal{O}(\varepsilon^7),\\
 \alpha &=& \frac{4 - \mu}{4} - \frac{1}{12}\varepsilon^2 + \frac{5}{3456 } \varepsilon^4 - \frac{169 249}{4 976 640} \varepsilon^6 + \mathcal{O}(\varepsilon^7).
\end{eqnarray*}

\begin{figure}[ht]
 \centering
 \includegraphics[scale=0.8]{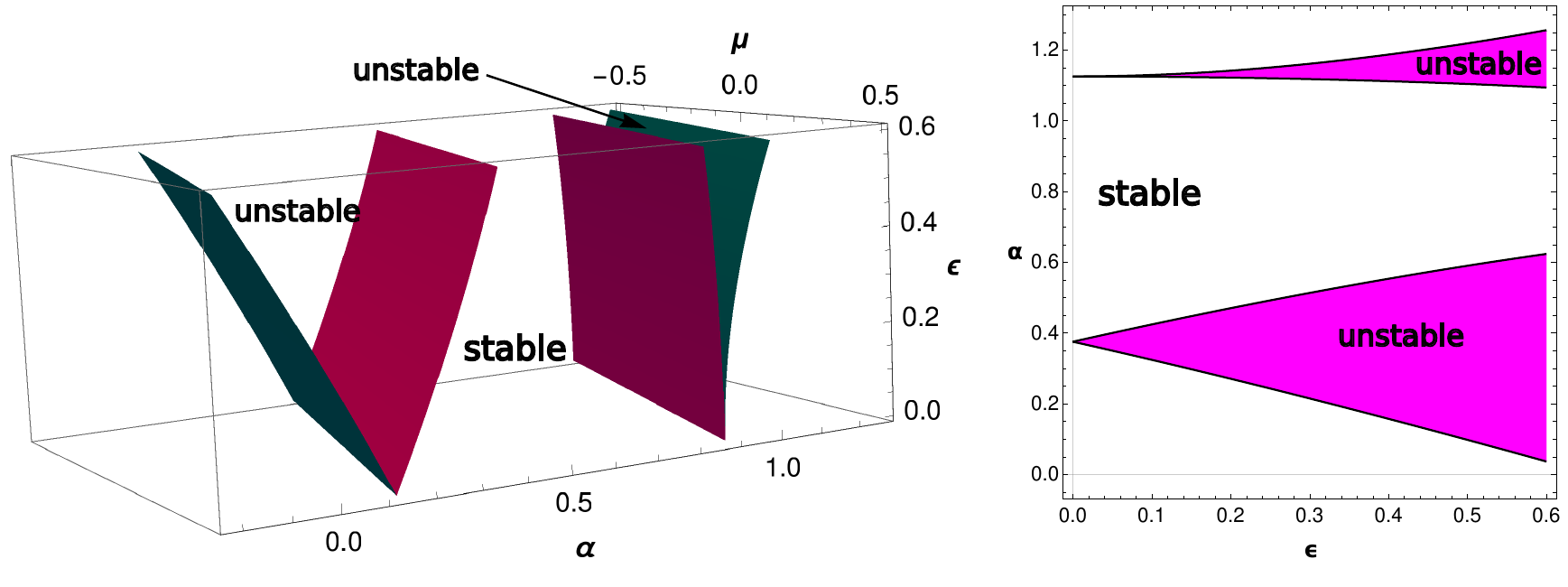}
  \caption{Boundary surfaces for $P_1$ with $N =1$ and, with $N = 2$; Planar section $\mu = -\frac{1}{2}.$}
  \label{fig3}
 \end{figure}

In Fig. \ref{fig3}, the picture on the left hand side shows the regions of stability and instability of the parameter space associated to the resonance $2\omega = 1 $ and $2 \omega = 2$ for the equilibrium $P_1$. The picture on the right hand side is a planar section of that regions with respect to the plane $\mu = -\frac{1}{2}$.

For $N = 3$ we obtain
\begin{eqnarray*}
 k_{20}^{(1)} &=& \frac{\alpha_1}{6},\\
 k_{20}^{(2)} &=& -\frac{1}{864}(9 + 24 \alpha_1 + 16 \alpha_1^2 - 144 \alpha_2),\\
 k_{20}^{(3)} &=& \frac{1}{15 552}(-108 + 153 \alpha_1 + 384 \alpha_1^2 + 64 \alpha_1^3 - 576 \alpha_1 \alpha_2 - 432 \alpha_2 + 2592 \alpha_3),
\end{eqnarray*}
and,
\begin{center}
$k_{02}^{(1)} = k_{20}^{(1)}$, \quad $k_{02}^{(2)} - k_{20}^{(2)} = \frac{\alpha_1}{18}$, \quad $k_{02}^{(3)} - k_{20}^{(3)} = -\frac{1}{648}(9 + 32 \alpha_1^2 - 36 \alpha_2)$. 
\end{center}
For $k_{20} = 0$ and $k_{02} = 0$, we obtain, respectively, the surfaces
$$\alpha = \frac{9 - \mu}{4} +\frac{1}{16}\varepsilon^2 \mp \frac{1}{32}\varepsilon^3 + \frac{13 }{5120}\varepsilon^4 \pm \frac{5 }{2048}\varepsilon^5 + \mathcal{O}(\varepsilon^6).$$

\begin{figure}[ht]
 \centering
 \includegraphics[scale= 0.9]{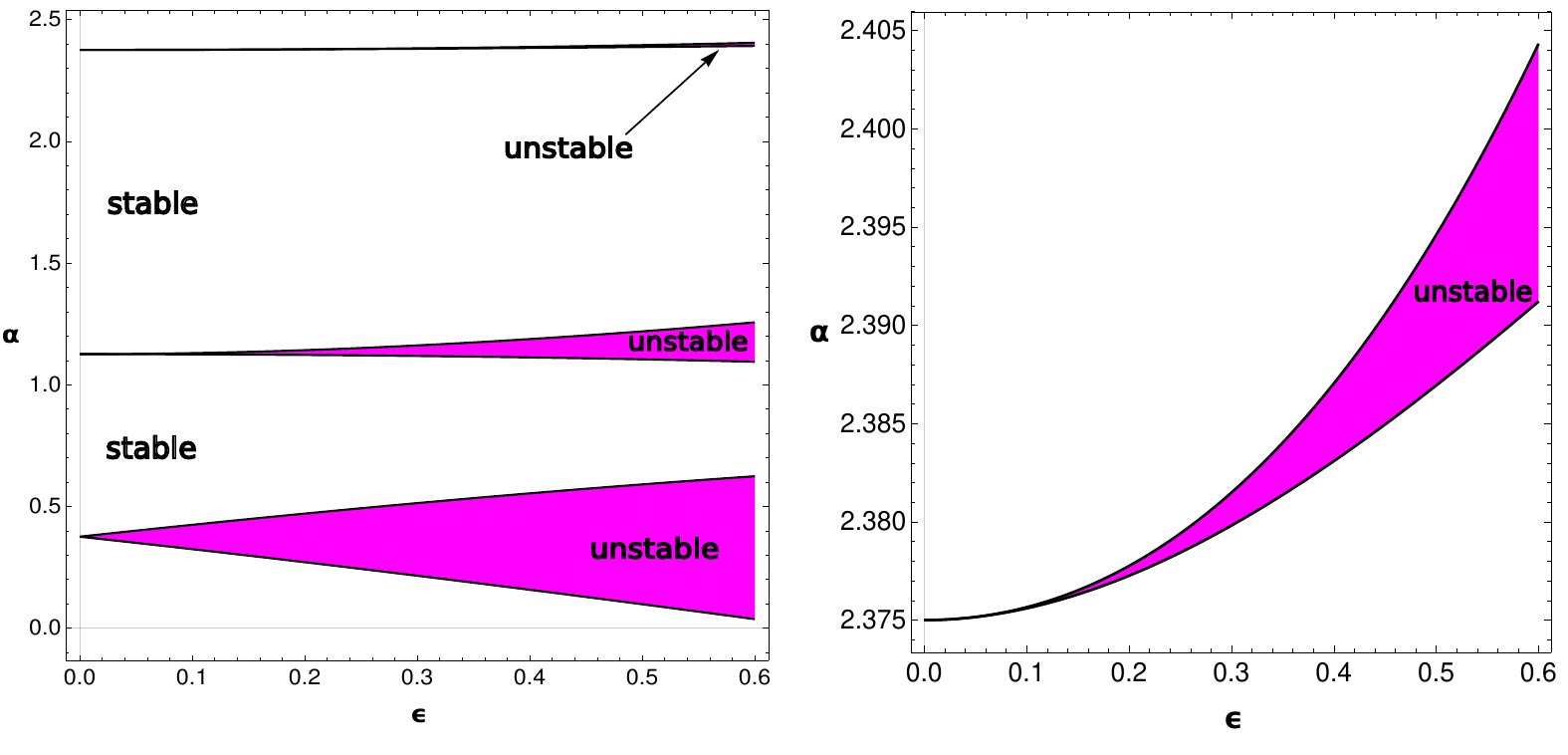}
  \caption{Planar section for $\mu = -\frac{1}{2}$; Expansion of the case $N = 3.$}
  \label{fig4}
 \end{figure}

We then follow the process for $N =  4, 5, 6, \cdots$, obtaining a decomposition of the parameter space $(\mu, \alpha, \varepsilon)$ interspersed by stability and instability regions for the equilibrium $P_1$. In Fig. \ref{fig4}, the left hand side picture shows the decomposition for $N =1, 2, 3 $ on the plane $\mu =-\frac{1}{2}$. The right hand side picture is an expansion of the case $N =3$.

\subsection{Boundary surfaces for the equilibrium P2}

When applying the Deprit-Hori Method on the Hamiltonian function \eqref{eq:hamiltonian4} for each resonance of the form $2 \omega = N$, $N = 1, 2, 3, \cdots$, with $\omega_0^2 = \frac{\mu}{4} - \alpha_0$, we obtain the following surfaces

For $N = 1$, the equations $k_{20} = 0$ and $k_{02} = 0$ provide, respectively the following surfaces
\[
\alpha = \frac{\mu - 1}{4} \mp \frac{1}{2} \varepsilon +\frac{1 }{8}\varepsilon^2  \pm \frac{1 }{32} \varepsilon^3 + \frac{1 }{384} \varepsilon^4 \mp \frac{11}{4608}  \varepsilon^5 + \mathcal{O}(\varepsilon^6).
\]

For $N = 2$, the surfaces obtained from the equations $k_{20} = 0$ and $k_{02} = 0$ are, respectively
\[
 \alpha = \frac{\mu - 4}{4}
\quad \text{ and } \quad \alpha = \frac{\mu - 4}{4}.
\]

\begin{figure}[ht]
 \centering
 \includegraphics[scale= 0.5]{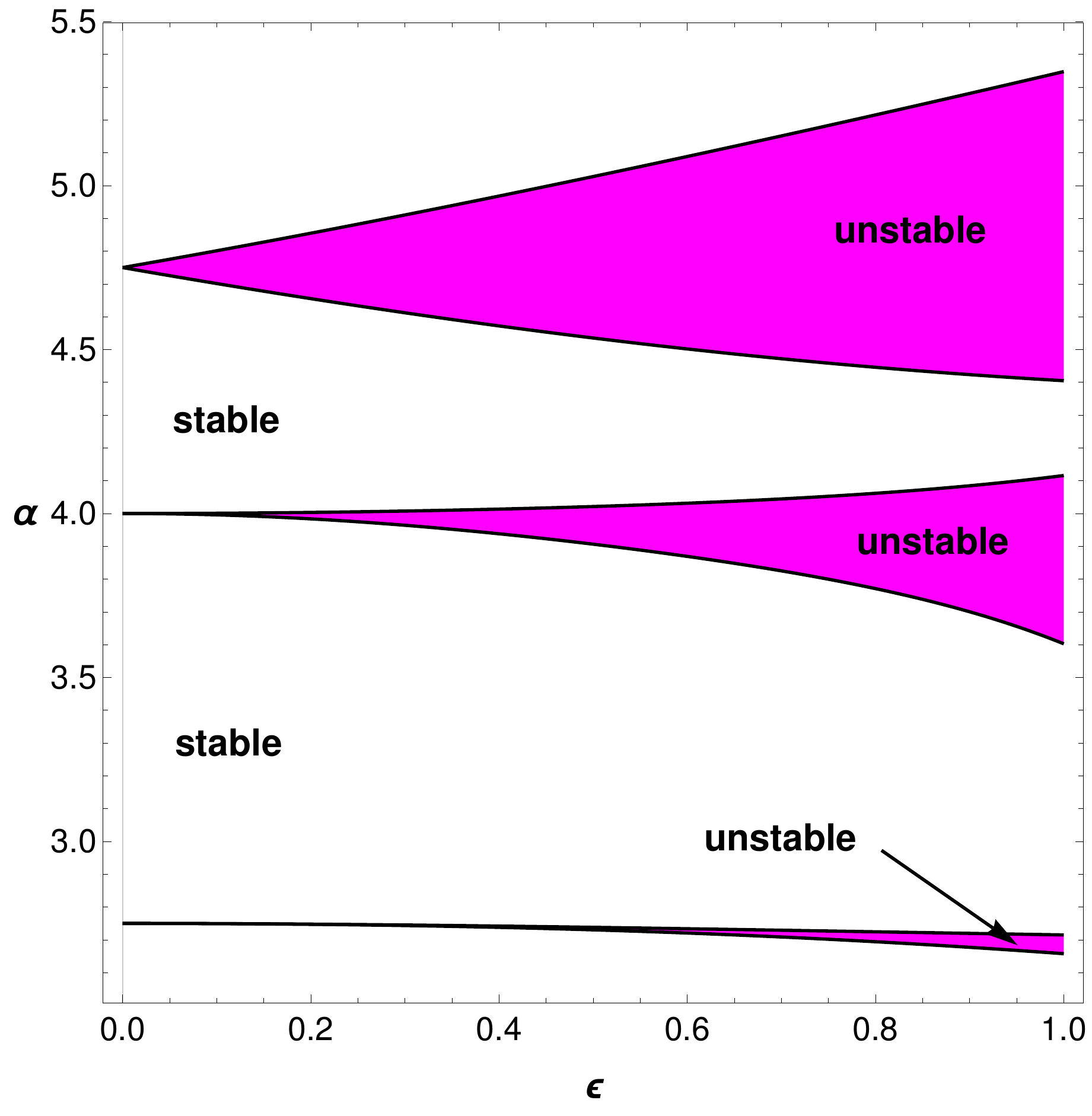}
  \caption{Planar section of the regions with respect to $\mu = 20.$}
  \label{fig5}
 \end{figure}

Finally, for $N = 3$, we have
\[
\alpha = \frac{\mu - 9}{4} - \frac{1}{16}\varepsilon^2 \mp\frac{1}{32}\varepsilon^3 - \frac{13 }{5120} \varepsilon^4 \pm \frac{5 }{2048} \varepsilon^5 + \mathcal{O}(\varepsilon^6).\]
Analogously to what we did for equilibrium $P_1$, following the process for $N =  4, 5, 6, \cdots$ we obtain a decomposition of the parameter space $(q, \alpha, \varepsilon)$ into stability and instability regions for the equilibrium $P_2$. Fig. \ref{fig5}  shows such decomposition for $N =1, 2, 3 $ restricted to the plane  $\mu = 20 $. 

\section{Conclusion}

In this work, we studied the parametric resonances of a model describing the dynamics of a mathematical pendulum with
support point oscillating vertically in a harmonic way under the influence of two electrically charged
lines with uniform charge distribution, equidistant from the pendulum support point.
By taking a Hamiltonian formulation, we determined the stability of the equilibria in the parameter space $(\mu, \alpha, \varepsilon)$. We proved that the equilibria $P_1 = (0,0)$ and $P_2= (\pi,0)$ are linearly stable for $\mu > 4 \alpha$ and $\mu > - 4\alpha$ respectively.
We normalized the quadratic part of the Hamiltonian function and then, we applied the Deprit-Hori Method in order to obtain the surfaces in the parameter space which separates the stability and the instability regions for the observed equilibria.
The surfaces were obtained as a graph of a function on $(\mu,0,\varepsilon)$ by determining the coefficients of its parametrizations until fifth order on $\varepsilon$, in terms of $\mu$.  The particular case $\mu = 0$ shows the boundary curves of the Mathieu equation, confirming the results published in \cite{Bardin_Markeev}.

\end{document}